\documentclass[10pt]{amsart}
\usepackage[cp1251]{inputenc}
\usepackage[english,russian]{babel}
\usepackage{amsmath}
\usepackage{amssymb}
\usepackage{amsfonts}
\usepackage{graphicx}

\begin{document}
\renewcommand{\refname}{References}
\renewcommand{\proofname}{Proof.}
\renewcommand{\figurename}{Fig.}

\sloppy

\title{Asymptotics of sums of regression residuals under multiple ordering of regressors}
\author{Mikhail Chebunin, Artyom Kovalevskii}
\thanks{E-mail: chebuninmikhail@gmail.com, 
Karlsruhe Institute of Technology, Institute of Stochastics, 76131 Karlsruhe,
Germany; Novosibirsk State University, Novosibirsk, Russia;
 E-mail: artyom.kovalevskii@gmail.com, Novosibirsk State Technical University, 
Novosibirsk State University, Novosibirsk, Russia.  
The work is supported by Mathematical Center in Akademgorodok 
under agreement No. 075-15-2019-1675 with the Ministry of Science and Higher Education of the Russian Federation.
}

\date{}
\maketitle 

\noindent{\sc Abstract}

We prove theorems about the Gaussian asymptotics of an empirical bridge built from linear model regressors with multiple regressor ordering. We study the testing of the hypothesis of a linear model for the components of a random vector: one of the components is a linear combination of the others up to an error that does not depend on the other components of the random vector. The results of observations of independent copies of a random vector are sequentially ordered in ascending order of several of its components. The result is a sequence of vectors of higher dimension, consisting of induced order statistics (concomitants) corresponding to different orderings. For this sequence of vectors, without the assumption of a linear model for the components, we prove a lemma of weak convergence of the distributions of an appropriately centered and normalized process to a centered Gaussian process with almost surely continuous trajectories. Assuming a linear relationship of the components, standard least squares estimates are used to compute regression residuals --- the differences between response values and the predicted ones by the linear model. We prove a theorem of weak convergence of the process of regression residuals under the necessary normalization to a centered Gaussian process. 

\medskip

\noindent{\bf Keywords:} concomitants, copula, weak convergence, regression residuals, empirical bridge.

\section{Introduction}

An extremely useful method for analyzing multivariate statistics
is the study of linear relationships between components. This analysis allows you to build
a linear prediction of one variable based on the others.
The very existence of dependencies is verified by calculating sample correlations and developing tests based on them.
This class of tests is the subject of correlation analysis.

The construction of models for the linear dependence of one variable (response) on other variables (regressors), the estimation of the parameters of the linear dependence and testing of their significance are the subject of
regression analysis. However, standard regression analysis methods do not include
methods of detecting that the proposed linear model is incorrect entirely. If the model is incorrect, then
it must either be completely discarded or substantially modified.

Methods for testing linear models, as a rule, use functionals as statistics
from random processes built according to the sequence of observations.
If the observations are ordered by one of the regressors, then such statistical tests are often called tests
discord detection. In the papers of Rao (1950), Page (1954), observations are ordered by time, and the alternative
hypothesis is that the distribution changes at some time (the change point). In this case, the distribution before the change
is assumed to be known, and the tests are focused on the fastest detection of the change.
Moustakides (1986) proved that the CUSUM procedure proposed by Page (1954)
is optimal in terms of Lorden (1971).
Shiryaev (1996) generalized this result to a continuous-time analogue of CUSUM.
Brodsky and Darkhovsky (2005) proved the asymptotic optimality of the adaptive CUSUM test for  compact
sets of unknown distribution parameters.

In situations where the distribution before the disorder is not described by a set of parameters from a compact set,
and also for more complex linear models, the process of sums of regression residuals is used, 
see Shorack and Wellner (1986).
MacNeill (1978) proposed such a test for time series, and Bishoff (1998) significantly relaxed the assumptions of
MacNeill.  An analysis of results in this direction can be found in Csorgo and Horv{\'a}th (1997, Chapters 2 and 3) and MacNeill et al. (2020). Kovalevskii and Shatalin (2015, 2016), Kovalevskii (2020) proposed tests for matching of regression
models using data ordering by one of the regressors. We offer a statistical test that uses
multiple ordering of data across multiple regressors.

The rest of the work is organized as follows. We prove a lemma extending
the functional central limit theorem by Davydov and Egorov (2000) for the multidimensional case (convergence
to a Gaussian field) in Section 2. This lemma is based on the general result of Ossiander (1987) and 
allows one to obtain a limit theorem for multiple ordering.
Section 3 contains this limit theorem for a general linear regression model, as well as an example of its application to the analysis of multivariate data. Section 4 contains the proofs of the Lemma and the Theorem.

\section{ Induced order statistics}

Let $\left({\bf X}_{i}, {\bf Y}_{i}\right)$, $i=1,2 \ldots,$ be the independent copies of a random vector 
$({\bf X}, {\bf Y})$ such that 
${\bf X}=(X^{(1)},\ldots, X^{(d_1)} )$ takes values in $[0,1]^{d_1}$, ${\bf Y}$ takes values in $\mathbb{R}^{d_2}$. 
The distribution function (copula) of ${\bf X}$ is 
$C({\bf u})={\bf P}({\bf X} \le {\bf u}) = {\bf P} \left({X}^{(1)}\le u^{(1)}, \ldots, {X}^{(d_1)}\le u^{(d_1)}\right)$, 
${\bf u}=(u^{(1)}, \ldots, u^{(d_1)}) \in [0,1]^{d_1}$.

We assume that there is copula density $c({\bf u})$, that is, 
\begin{equation}\label{copula_density}
C({\bf u})=\int_{{\bf v}\le {\bf u}} c({\bf v}) d {\bf v}, \ \ {\bf v}\in [0,1]^{d_1}.
\end{equation}

Denote $X^{(k)}_{n, 1} \leq X^{(k)}_{n, 2} \leq \cdots \leq
X^{(k)}_{n, n}$, $1\le k\le d_1$, the order statistics of the $k$-th column of matrix $X$, and 
${\bf Y}^{(k)}_{n, 1}, {\bf Y}^{(k)}_{n, 2}, \ldots, {\bf Y}^{(k)}_{n, n}$ the corresponding values of the vectors 
${\bf Y}_i$. The random vectors $\left({\bf Y}^{(k)}_{n, i}, i \leq n\right)$ are called induced order statistics (concomitants).

We study the asymptotic behavior of random field
$$
\begin{aligned}
{\bf Q}_{n}({\bf u}) &= \sum_{j=1}^{n} {\bf Y}_{ j} \mathbf{1}\left({\bf X}_{ j}\le {\bf u}\right)
= \sum_{j=1}^{n} {\bf Y}_{ j} \mathbf{1}\left({X}_{ j}^{(1)}\le u^{(1)}, \ldots, {X}_{ j}^{(d_1)}\le u^{(d_1)}\right), \ {\bf u}\in [0,1]^{d_1}.
\end{aligned}
$$

Using the asymptotics of ${\bf Q}_{n}({\bf u})$, we study the asymptotics of $d_1\times d_2$-dimensional process 
of sums of induced order statistics under different orderings
$$
\begin{aligned}
{\bf Z}_{n}(t) &= \left(
\sum_{j=1}^{[n t]} {\bf Y}^{(1)}_{n, j},
\sum_{j=1}^{[n t]} {\bf Y}^{(2)}_{n, j},
 \dots,
 \sum_{j=1}^{[n t]} {\bf Y}^{(d_1)}_{n, j} 
 \right), \ t \in [0, 1].
\end{aligned}
$$

Let ${\bf m}({\bf u})={\bf E}({\bf Y} \mid {\bf X}={\bf u})$, ${\bf u} \in  [0,1]^{d_1}$, and ${\bf f}({\bf u})=\int_{{\bf 0}}^{\bf u} {\bf  m}({\bf v}) c({\bf v}) d {\bf v}$.

 Let
$$
\sigma^{2}({\bf u})={\bf E}\left\{({\bf Y}-{\bf m}({\bf X}))^{T}({\bf Y}-{\bf m}({\bf X})) \mid {\bf X}={\bf u}\right\}
$$
be the conditional covariance matrix of ${\bf Y}$ and $\sigma({\bf u})$ be the positive definite matrix such that $\sigma({\bf u})^{T} \sigma({\bf u})=\sigma^{2}({\bf u})$.

Let ${\bf e}_{k,t}=(1,\ldots,1,t,1,\ldots,1)$ the vector in $[0,1]^{d_1}$ with $k$-th coordinate being $t$ and other
coordinates being 1.

All our limit fields and processes being continuous, so we  use the uniform metric. The symbol $\Rightarrow$ denotes the weak convergence of random variables or weak convergence of random fields or stochastic processes in the uniform topology. Let $\|\cdot\|$ denote Euclidean norm in the corresponding space.

The following Lemma 1 generalizes the result of the first part of Theorem 2.1(1) by Davydov and Egorov (2000) to random fields.

{\bf Lemma 1} If ${\bf E} \|{\bf Y}\|^2 < \infty$ then
$\widetilde{\bf Q}_n=\frac{{\bf Q}_n - {\bf f}}{\sqrt{n}} \Rightarrow {\bf Q}$, a centered Gaussian field with covariance 
\[
K({\bf u}_1, {\bf u}_2) = {\bf E} {\bf Q}^T({\bf u}_1) {\bf Q}({\bf u}_2) = 
\int_{\bf 0}^{\min({\bf u}_1,{\bf u}_2)} \sigma^2({\bf v}) c({\bf v}) d {\bf v} 
\]
\[
+ 
\int_{\bf 0}^{\min({\bf u}_1,{\bf u}_2)} {\bf m}^T({\bf v}) {\bf m}({\bf v}) c({\bf v}) d {\bf v}
- 
\int_{\bf 0}^{{\bf u}_1} {\bf m}^T({\bf v}) c({\bf v}) d {\bf v}
\int_{\bf 0}^{{\bf u}_2} {\bf m}({\bf v}) c({\bf v}) d {\bf v},
\]
${\bf u}_1$, ${\bf u}_2 \in [0,1]^{d_1}$;

Lemma 2 generalizes the result of  Theorem 2.1(2) by Davydov and Egorov (2000) to 
multiple ordering but under additional assumption ${\bf m}\equiv {\bf 0}$.

{\bf Lemma 2} If ${\bf E} \|{\bf Y}\|^2 < \infty$, ${\bf m}\equiv {\bf 0}$ then $\widetilde{\bf Z}_n=\frac{{\bf Z}_n}{\sqrt{n}} \Rightarrow {\bf Z}$, a centered Gaussian $(d_1\times d_2)$-dimensional process with
covariance matrix function  ${\bf E}{\bf Z}^T(t_1){\bf Z}(t_2)=(K({\bf e}_{k_1,t_1}, {\bf e}_{k_2,t_2}))_{k_1,k_2=1}^{d_1}$, 
\[
K({\bf e}_{k_1,t_1}, {\bf e}_{k_2,t_2}) = {\bf E} 	{\bf Q}^T({\bf e}_{k_1,t_1}){\bf Q}({\bf e}_{k_2,t_2}) = 
\int_{\bf 0}^{\min({\bf e}_{k_1,t_1},{\bf e}_{k_2,t_2})} \sigma^2({\bf v}) c({\bf v}) d {\bf v}. 
\]

\vspace{5 mm}

\section{Main result}

Let $({\bf X}_i, \xi_i, \eta_i) =(X_{i1},\ldots,X_{id_1},\xi_{i1},\ldots, \xi_{i,d_2-1}, \eta_i)$ be independent and identically distributed 
random vector rows, $i=1,\ldots,n$. 
All components of a raw can be dependent and $X_{i1},\ldots,X_{id_1}$ have copula (so their marginal distributions are 
uniform on [0, 1]) and (\ref{copula_density}) is true.

Rows $({\bf X}_i, \xi_i, \eta_i)$ form matrix $(X,\xi, \eta)$.

We assume a linear regression hypothesis $H_0$: 
\begin{equation}\label{regression}
\eta_i=\xi_i \theta+\varepsilon_i =
\sum_{j=1}^{d_2-1} \xi_{ij}\theta_j+\varepsilon_i,
\end{equation}
$\{\varepsilon_i\}_{i=1}^n$ and $\{({\bf X}_i, \xi_i)\}_{i=1}^n$ are independent, ${\bf E}\, \varepsilon_1 =0$, 
${\bf Var} \, \varepsilon_1>0$.

Vector $\theta=(\theta_1,\ldots,\theta_{d_2-1})^T$ and constant ${\bf Var} \, \varepsilon_1$ are unknown. 
We consider $d_1$ orderings of rows of the matrix $ (X,\xi, \eta) $ in acsending order of 
columns  of $X$.

The result of $d_1$ orderings is a sequence of $d_1$ matrices $(X^{(j)},\xi^{(j)},\eta^{(j)})$ 
with rows $({\bf X}_i^{(j)},{\xi}_i^{(j)}, \eta_i^{(j)})=(X_{i1}^{(j)},\ldots, X_{id_1}^{(j)}, 
\xi_{i1}^{(j)}, \ldots, \xi_{i,d_2-1}^{(j)}, \eta_i^{(j)})$, $j=1,\ldots,d_1$. 
 
\bigskip
 
Let $\widehat{\theta}$ be LSE:
\[
\widehat{{\theta}}=({\xi}^T\xi)^{-1}{\xi}^T{\eta}.
\]

It does not depend on the order of rows.

Let
$h^{(j)}(x)={\bf E}\{\xi_{1} | X_{1j}=x\}$ be conditional expectations,
$L^{(j)}(x)=\int\limits_{0}^{x}h^{(j)} (s)\,ds$  be {\it induced 
theoretical generalised Lorentz curves} (see Davydov and Egorov (2000)),
\[
b^2_j(x)={\bf E} \left((\xi_1-h^{(j)}(x))^T (\xi_1-h^{(j)}(x)) \ |\ X_{1j}=x\right)
\]
be matrices of conditional covariances.

Let $G={\bf E}\xi_{1}^T\xi_{1}$.
Then 
\[
\int_0^1\left(b_j^2({ x})+(h^{(j)}({ x}))^T h^{(j)}({ x})\right) \, d{ x}=G
\]
 for any $j=1,\ldots,d_2-1$.

Let $\widehat{\varepsilon}_i^{(j)}=\eta_i^{(j)} - \xi_i^{(j)} \widehat{\theta}$ be regression residuals,
$
\widehat{\Delta}_k^{(j)}=\sum\limits_{i=1}^{k}\widehat{\varepsilon}_i^{(j)}$ be its partial sums, 
$\widehat{\Delta}_0^{(j)}=0$.

Let $\widehat{Z}_n^{(j)}=\{\widehat{Z}_n^{(j)}(t), \, 0 \le t \le 1\}$ be a piecewise linear random function with nodes
\[
\left(\frac{k}{n}, \
\frac{\widehat{\Delta}_k^{(j)}}{\sqrt{n {\bf Var} \varepsilon_1}}\right), \ \ k=0, \ 1, \ldots, \ n.
\]

 From Theorem 1 (Kovalevskii, 2020) we have

{\bf Theorem 1} {\em If matrix $G$ exists and is non-degenerate and $H_0$ is true
then $\widehat{Z}_n^{(j)}\Longrightarrow~\widehat{Z}^{(j)}$ for any $j=1,\ldots,d_1$. Here $\widehat{Z}^{(j)}$ 
is a centered Gaussian process with  continuous a.s. sample paths and covariance function }
\[
\widehat{K}_{jj}(s,t)=\min(s,t)- L^{(j)}(s) G^{-1} (L^{(j)}(t))^T, \
~s, 
t\in[0,1].
\]

We prove that the $d_1$-dimensional process $\widehat{Z}_n=(\widehat{Z}_n^{(j)}, j=1,\ldots,d_1)$ has a Gaussian limit.

{\bf Theorem 2} {\em If matrix $G$ exists and is non-degenerate and $H_0$ is true
then $\widehat{Z}_n\Longrightarrow~\widehat{Z}$. Here $\widehat{Z}$ is a centered $d_1$-dimensional Gaussian process with continuous a.s. sample paths and covariance matrix function } 
$\widehat{K}(s,t)=\left(\widehat{K}_{ij}(s,t)\right)_{i,j=1}^{d_1}$,
\[
\widehat{K}_{ij}(s,t)={\bf P}(X_{1i}
\le s, X_{1j} \le t) - L^{(i)}(s) G^{-1} (L^{(j)}(t))^T, \
~s, 
t\in[0,1].
\]

\section{Proofs}

{\it Proof of Lemma 1}

For simplicity, we consider the case $ d_1 = 2 $ since the construction of the proof given below can be easily extended to the case $ d_1> 2 $. Now let $d_2=1$, we will generalize it to $d_2\ge 1$ using the Cramer-Wold theorem.

Thus, we consider a random field
\[
Q_{n}({\bf u}) = \sum_{j=1}^{n}  Y_{ j} \mathbf{1}\left({X}_{ j}^{(1)}\le u^{(1)}, {X}_{ j}^{(2)}\le u^{(2)}\right), \ {\bf u}\in [0,1]^{2}.
\]

Let us define the partition of the unit square   $[0,1]^2$ into $N^2$ parts as follows.
Let $u^{(1)}_{0}=0<u^{(1)}_{1}<u^{(1)}_{2}<\ldots<u^{(1)}_{N}=1$ be a partition of the interval [0,1], such that
\[
\int_{(u^{(1)}_{i-1},0)}^{(u^{(1)}_i,1)}\left({\bf E}\left(Y^{2} \mid {\bf X}={\bf v}\right)+{\bf E} Y^{2}\right) c({\bf v}) d {\bf v}=2 {\bf E} Y^{2}/N, \ i=1,2, \ldots, N,
\]
and for any fixed $i=1,2, \ldots, N$ let $u^{(2)}_{i,0}=0<u^{(2)}_{i,1}<u^{(2)}_{i,2}<\ldots<u^{(2)}_{i,N}=1$ be another partition of the interval [0,1] (see Pic.1) such that
\[
\int_{(u^{(1)}_{i-1},u^{(2)}_{i,j-1})}^{(u^{(1)}_i,u^{(2)}_{i,j})}\left({\bf E}\left(Y^{2} \mid {\bf X}={\bf v}\right)+{\bf E} Y^{2}\right) c({\bf v}) d {\bf v}=2 {\bf E} Y^{2}/N^2, \ j=1,2, \ldots, N.
\]

\begin{center}
{\unitlength = 2 pt
\begin{picture}(150.00,140.00)

\thicklines

\put(20.00,20.00){\vector(1,0){130.00}}
\put(20.00,20.00){\vector(0,1){110.00}}
\put(12.00,130.00){$u^{(2)}$}
\put(150.00,14.00){$u^{(1)}$}
\put(20.00,120.00){\line(1,0){100.00}}
\put(120.00,20.00){\line(0,1){100.00}}
\put(50.00,20.00){\line(0,1){100.00}}
\put(100.00,20.00){\line(0,1){100.00}}
\put(20.00,40.00){\line(1,0){30.00}}
\put(20.00,70.00){\line(1,0){30.00}}
\put(20.00,40.00){\line(1,0){30.00}}
\put(50.00,80.00){\line(1,0){50.00}}
\put(50.00,100.00){\line(1,0){50.00}}
\put(100.00,30.00){\line(1,0){20.00}}
\put(100.00,60.00){\line(1,0){20.00}}
\put(16.00,14.00){$0$}
\put(50.00,14.00){$u^{(1)}_1$}
\put(100.00,14.00){$u^{(1)}_2$}
\put(120.00,14.00){$u^{(1)}_3=1$}
\put(11.00,40.00){$u^{(2)}_{1,1}$}
\put(11.00,70.00){$u^{(2)}_{1,2}$}
\put(2.00,120.00){$u^{(2)}_{1,3}=1$}

\end{picture}
}

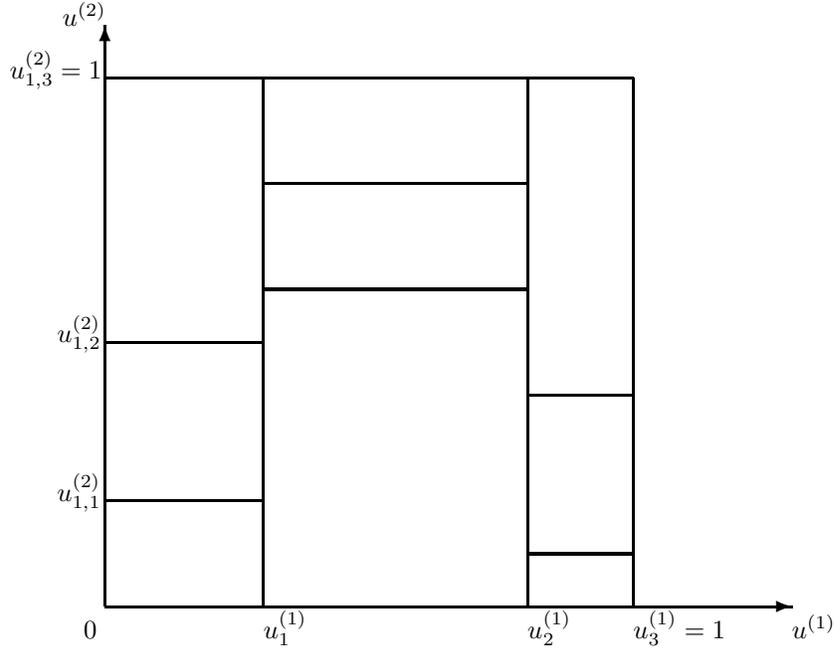
Pic. 1. An example of the partition of $[0,1]^2$ for $N=3$.

\end{center}

So we have $N+1$ points $u^{(1)}_{i}$ in the first coordinate and not greater then $(N-1)^2+1$ different points $u^{(2)}_{i,j}$ in the second coordinate. For any ${\bf u}=(u^{(1)},u^{(2)})\in[0,1]^2$, there are indexes $i^{(1)},i^{(2)}_1,i^{(2)}_2,j^{(2)}_1,j^{(2)}_2\in\{0,1,\dots,N\}$ such that
$u^{(1)}_{i^{(1)}-1}\le u^{(1)}\le u^{(1)}_{i^{(1)}}$ and
\[
u^{(2)}_{i^{(2)}_1,j^{(2)}_1}=\max_{i,j}\{ u^{(2)}_{i,j}\le u^{(2)}\},
\]
\[
u^{(2)}_{i^{(2)}_2,j^{(2)}_2}=\min_{i,j}\{ u^{(2)}_{i,j}\ge u^{(2)}\},
\]
so
$u^{(2)}_{i^{(2)}_1,j^{(2)}_1}\le u^{(2)}\le u^{(2)}_{i^{(2)}_2,j^{(2)}_2}$.

 Denote ${\bf u}^l=(u^{(1)}_{i^{(1)}-1}, u^{(2)}_{i^{(2)}_1,j^{(2)}_1})$ and 
${\bf u}^u=( u^{(1)}_{i^{(1)}},u^{(2)}_{i^{(2)}_2,j^{(2)}_2})$ (Pic. 2).

\begin{center}
{\unitlength = 2 pt
\begin{picture}(150.00,140.00)

{\thicklines

\put(20.00,20.00){\vector(1,0){130.00}}
\put(20.00,20.00){\vector(0,1){110.00}}
\put(12.00,130.00){$u^{(2)}$}
\put(150.00,14.00){$u^{(1)}$}
\put(20.00,120.00){\line(1,0){100.00}}
\put(120.00,20.00){\line(0,1){100.00}}
\put(50.00,20.00){\line(0,1){100.00}}
\put(100.00,20.00){\line(0,1){100.00}}
\put(20.00,40.00){\line(1,0){30.00}}
\put(20.00,70.00){\line(1,0){30.00}}
\put(20.00,40.00){\line(1,0){30.00}}
\put(50.00,80.00){\line(1,0){50.00}}
\put(50.00,100.00){\line(1,0){50.00}}
\put(100.00,30.00){\line(1,0){20.00}}
\put(100.00,60.00){\line(1,0){20.00}}
\put(16.00,14.00){$0$}
\put(50.00,14.00){$u^{(1)}_1$}
\put(100.00,14.00){$u^{(1)}_2$}
\put(120.00,14.00){$u^{(1)}_3=1$}
\put(11.00,40.00){$u^{(2)}_{1,1}$}
\put(11.00,70.00){$u^{(2)}_{1,2}$}
\put(2.00,120.00){$u^{(2)}_{1,3}=1$}
}

\put(20.00,90.00){\line(1,0){90.00}}
\put(110.00,20.00){\line(0,1){70.00}}

\put(110,90){$\bf u$}
\put(100,80){${\bf u}^{\it l}$}
\put(120,100){${\bf u}^{\it u}$}

\end{picture}
}

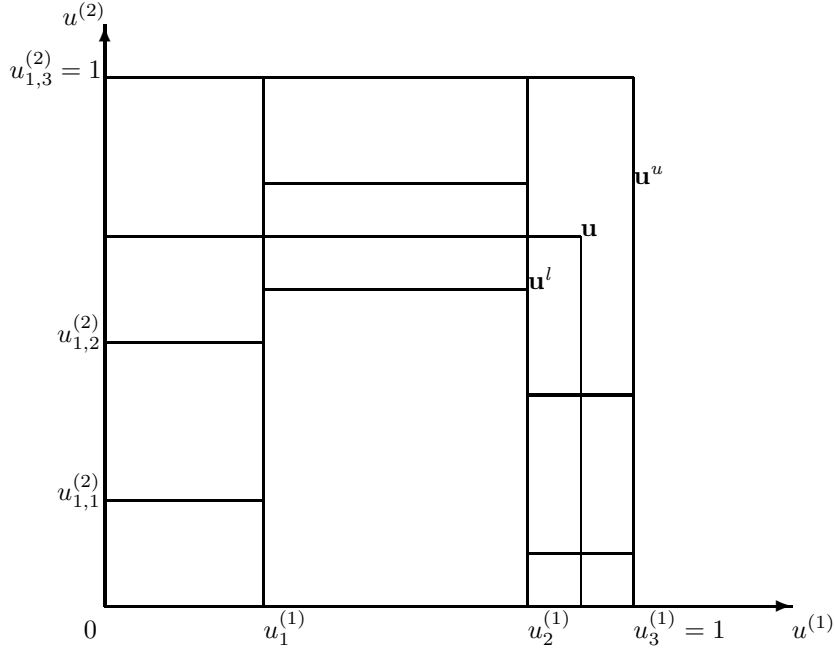
Pic. 2. An example of upper and lower points for $N=3$.

\end{center}

Define the metric entropy with bracketing for the special separable pseudometric space $(S, \rho),$ where
$$
\begin{array}{l}
S=\left\{h_{\bf u}({\bf x}, y)=y \mathbf{1}({\bf x}\le {\bf u}), \quad {\bf u} \in[0,1]^2\right\} \\
\rho^{2}\left(h_{\bf u_1}, h_{\bf u_2}\right)={\bf E}\left(h_{\bf u_1}({\bf X}, Y)-h_{\bf u_2}({\bf X}, Y)\right)^{2}
\end{array}
$$
Let $S(\delta)=\left\{s_{1}, s_{2}, \ldots, s_M\right\} \subseteq S$ be such that for some random variables $f^{l}\left(s_{i}\right)$ and $f^{u}\left(s_{i}\right), i \leq M,$ the following conditions are valid. For any $s \in S$ there exists $s_{i} \in S(\delta)$ such that
$$
\begin{aligned}
\rho\left(s, s_{i}\right) & \leq \delta, \\
f^{l}\left(s_{i}\right) \leq f(s,{\bf X}, Y) & \leq f^{u}\left(s_{i}\right) \ \ a . s ., \\
\rho\left(f^{u}(s_{i}), f^{l}(s_{i})\right) & \leq \delta, \ \ i \leq n.
\end{aligned}
$$
Then $H^{B}(\delta, S, \rho)=\min \{ M: S(\delta) \subseteq S\}$ is called the metric entropy with bracketing.

According to Ossiander's (1987) theorem, to prove Lemma for $\widetilde{ Q}_n$ we must show that
$$
\int_{0}^{1}\left(H^{B}(t, S, \rho)\right)^{1 / 2} d t<\infty.
$$
To prove it we show that
$$
H^{B}(\delta, S, \rho) \leq C \log \left(\frac{1}{\delta}\right).
$$
First notice that $f\left(h_{\bf u}, {\bf X}, Y\right)=h_{\bf u}({\bf X}, Y)-E h_{\bf u}({\bf X}, Y), {\bf u} \in [0,1]^2$ is a separable random process satisfying the conditions $(2.1)-(2.3)$ of Ossiander's work, moreover $\widetilde{ Q}_n=n^{-1 / 2} \sum_{i=1}^{n} f\left(h_{\bf u},{\bf X}_{i}, Y_{i}\right)$. Then write $f(h_{\bf u},{\bf x}, y)$ as follows
\[
f(h_{\bf u},{\bf x}, y)=y_{+} \mathbf{1}({\bf x}\le{\bf u})-y_{-} \mathbf{1}({\bf x}\le{\bf u})-\int_{{\bf 0}}^{{\bf u}} m_{+}({\bf v})c({\bf v}) d {\bf v}+\int_{{\bf 0}}^{{\bf u}} m_{-}({\bf v})c({\bf v}) d {\bf v}
\]

Let
\[
f^{u}=Y_{+} \mathbf{1}({\bf X}\le{\bf u}^u)-Y_{-} \mathbf{1}({\bf X}\le{\bf u}^l)
-\int_{{\bf 0}}^{{\bf u}^l} m_{+}({\bf v}) c({\bf v}) d {\bf v}+\int_{{\bf 0}}^{{\bf u}^u} m_{-}({\bf v}) c({\bf v}) d {\bf v},
\]
\[
f^{l}=Y_{+} \mathbf{1}({\bf X}\le{\bf u}^l)-Y_{-} \mathbf{1}({\bf X}\le{\bf u}^u)
-\int_{{\bf 0}}^{{\bf u}^u} m_{+}({\bf v}) c({\bf v}) d {\bf v}+\int_{{\bf 0}}^{{\bf u}^l} m_{-}({\bf v}) c({\bf v}) d {\bf v}.
\]

Then
\[
f^{u}-f^{l}=|Y|\mathbf{1}({\bf u}^l\le{\bf X}\le{\bf u}^u)+\int_{{\bf u}^l}^{{\bf u}^u} |m({\bf v})| c({\bf v})d {\bf v}.
\]

We use the Cauchy-Bunyakovsky and Jensen inequalities, as well as simple algebraic inequalities, and obtain that
(see Davydov and Egorov (2000) for details)
\[
\left(E\left(f^{u}-f^{l}\right)^{2}\right)^{1 / 2} \leq\left(\int_{{\bf u}^l}^{{\bf u}^u} E\left(Y^{2} \mid {\bf X}={\bf v}\right) c({\bf v})d {\bf v}\right)^{1 / 2}
\]
\[
+\left(E Y^{2}\int_{{\bf u}^l}^{{\bf u}^u} c({\bf v}) d {\bf v}\right)^{1 / 2}  \leq 2\left(\int_{{\bf u}^l}^{{\bf u}^u}\left(E\left(Y^{2} \mid {\bf X}={\bf v}\right)+E Y^{2}\right) c({\bf v}) d {\bf v}\right)^{1 / 2}
\]
\[
\le 2\sqrt{2N\cdot 2 E Y^{2}/N^2}=4\sqrt{E Y^{2}/N}=\delta.
\]

In the last inequality, we used the fact that the region of integration is included in $ 2N $ rectangles from the constructed partition of the unit square.

Summing up these equalities we get
$$
M< (N+1)^3 \leq \left(\left[\frac{16 E Y^{2}}{\delta^{2}}\right]+2\right)^3.
$$

Hence by Ossiander's theorem
$$
\widetilde{Q}_n \Rightarrow Q
$$
where $ Q$ is the Gaussian field, $E Q({\bf u})=0, {\bf E}(Q({\bf u}_1),Q({\bf u}_2))=K({\bf u}_1, {\bf u}_2)$.
Elementary calculations show that 
\[
K({\bf u}_1, {\bf u}_2)={\bf E} (Y^2 \mathbf{1}({\bf X}\le{\bf u}_1,{\bf X}\le{\bf u}_2)) - 
\int_{\bf 0}^{{\bf u}_1}  m({\bf v}) c({\bf v}) d {\bf v}
\int_{\bf 0}^{{\bf u}_2}  m({\bf v}) c({\bf v}) d {\bf v}.
\]

This construction of the proof can be easily extended to the case $ d_1> 2 $ by splitting the corresponding integral into pieces of size  $2{\bf E}Y^2/N^{d_1}$.

Note that since the trajectories of the limit Gaussian field are continuous, convergence in the Skorokhod metric is equivalent to convergence in the uniform metric. Therefore, the convergence in the uniform metric of all coordinate fields implies their relative compactness in the Skorokhod topology, and hence the relative
compactness of the initial random field. We have proved the convergence in the uniform metric
  for $ Q_{n} ({\bf u}) $ for $ d_2 = 1 $, and according to the Cramer-Wold theorem
since $ {\bf Q}_{n} ({\bf u}) $ is linear, we obtain the convergence of finite-dimensional distributions
for any $ d_2 \ge 1 $.

The proof is complete.

{\it Proof of Lemma 2}

 Note that
\[
{\bf Z}_{n}(t) = \left(Q_n\left({\bf e}_{1,X^{(1)}_{n,[nt]+1}}\right),  Q_n\left({\bf e}_{2,X^{(2)}_{n,[nt]+1}}\right), \dots,  Q_n\left({\bf e}_{d_1,X^{(d_1)}_{n,[nt]+1}}\right) \right), \ t \in [0, 1].
\]

Hence, due to Lemma 4.1 by Davydov and Egorov (2000), for all $i=1,\dots, d_1$
\[
\sup_{0<t<1}\left|\widetilde{Q}_n\left({\bf e}_{i,X^{(i)}_{n,[nt]+1}}\right)-\widetilde{Q}_n\left({\bf e}_{i,t}\right)\right| \stackrel{p}{\rightarrow} 0.
\]

We assume $m\equiv  0$, so
\[
\frac{1}{\sqrt{n}}\sup_{0<t<1}\left|Q_n\left({\bf e}_{i,X^{(i)}_{n,[nt]+1}}\right)-Q_n({\bf e}_{i,t})\right| \stackrel{p}{\rightarrow} 0.
\]

Therefore the limiting process for $\widetilde{Z}_n$ is the same as for  
\[
\frac{1}{\sqrt{n}}\left(Q_n({\bf e}_{1,t}),Q_n({\bf e}_{2,t}),\dots,Q_n({\bf e}_{d_1,t})\right).
\]

We have 
\[
K({\bf e}_{k_1,t}, {\bf e}_{k_2,t})={\bf E} (Y^2 \mathbf{1}({\bf X}\le{\bf e}_{k_1,t},{\bf X}\le{\bf e}_{k_2,t})).
\]

The proof is complete.

\bigskip

{\it Proof of Theorem 2}

Let $\varepsilon_i^{(j)}=\eta_i^{(j)}-\xi_i^{(j)}\theta$ be regression mistakes.
From (\ref{regression}) we have $\{\varepsilon_i^{(j)}\}_{i=1}^n$ are i.i.d. with $\varepsilon_1$ and independent with $\{(X^{(j)}_i, 
\xi_i^{(j)})\}$ for any $j=1,\ldots,d_1$.
 
Let $\varepsilon^{(j)}=(\varepsilon_1^{(j)},\ldots,\varepsilon_n^{(j)})^T$.
Note that 
\[
\widehat{\Delta}^{(j)}_k=\sum_{i=1}^k (\eta_i^{(j)} - \xi_i^{(j)} \widehat{\theta})
=
\sum_{i=1}^k (\xi_i^{(j)}(\theta- \widehat{\theta})+\varepsilon_i^{(j)})
\]
\[
=
\sum_{i=1}^k (\xi_i^{(j)}(\theta-({\xi}^T\xi)^{-1}{\xi}^T{\eta})+\varepsilon_i^{(j)})
\]
\[
=
\sum_{i=1}^k \left(\xi_i^{(j)}\left(\theta- ({\xi}^T\xi)^{-1}\left({\xi}^{(j)}\right)^T({\xi}^{(j)}\theta
+{\varepsilon}^{(j)})\right)+\varepsilon_i^{(j)}\right)
\]
\[
=
\sum_{i=1}^k \left(\varepsilon_i^{(j)}-\xi_i^{(j)}({\xi}^T\xi)^{-1}\left({\xi}^{(j)}\right)^T{\varepsilon}^{(j)}\right).
\]

Note that 
\[
\left\{
\sum_{i=1}^{[nt]} \xi_i^{(j)}/n, \ t \in [0,1] 
\right\} \to L^{(j)}
\] 
a.s. uniformely,
and ${\xi}^T{\xi}/n \to G$ a.s.

So we study process
\[
\left\{
\sum_{i=1}^{[nt]} \varepsilon_i^{(j)}-L^{(j)}(t) G^{-1} \left(\xi^{(j)}\right)^T{\varepsilon}^{(j)}, \ \
 \ \  1 \le j \le d_1, \ \ t\in [0,1]
\right\}.
\]

This process is a bounded linear functional of a $d_1\times d_2$-dimensional process
\[
\left\{
\sum_{i=1}^{[nt]} ({\xi}_i^{(j)}\varepsilon_i^{(j)}, \ \varepsilon_i^{(j)}), \ \
 \ \  1 \le j \le d_1, \ \ t\in [0,1]
\right\}.
\]

This is a process from Lemma 2 with ${\bf m}\equiv 0$. So we have convergence to a Gaussian process and calculate covariances using Lemma 2.

The proof is complete.

\section{Discussion}

We now describe the application of this result to testing the hypothesis of linear dependence.
Let $d_2-1=d_1$. Let $(\xi_i, \eta_i) =(\xi_{i1},\ldots, \xi_{i,d_1}, \eta_i)$ be independent and identically distributed 
random vector rows, $i=1,\ldots,n$. In addition, we  assume that the column $ \xi_{i, d_1} $ consists of ones:\begin{equation}\label{ones}
\xi_{i,d_1}\equiv 1.
\end{equation}

We want to test the linear dependence (\ref{regression}). To do this we estimate the parameters $ \theta $ and
$ {\bf Var} \, \varepsilon_1 $, sort the data in ascending order of each of the first $ d_1-1 $ columns of the regressor and calculate processes of the sums of regression residuals. We apply Theorem 2. We use the quantile functions 
$F^{-1}_{\xi_{1j}}$ for this.

We assume that $\xi_{ij}=F^{-1}_{\xi_{1j}}(X_{ij})$, $i=1,\ldots,n$, $j=1,\ldots,d_1-1$. 
If the matrix  $G$ exists and is non-degenerate then
under the true hypothesis $ H_0 $ we are in the conditions of Theorem 2. 

From (\ref{ones}) we have $\widehat{Z}_n(1)={\bf 0}$. So we can use a statistics of omega squared type and calculate its
limiting distribution by lines of Chakrabarty et al. (2020):
\[
\omega^2_n = \sum_{j=1}^{d_1-1} \int_0^1 \left(\widehat{Z}_n^{(j)}(t)\right)^2 dt
\Rightarrow
\omega^2 = \sum_{j=1}^{d_1-1} \int_0^1 \left(\widehat{Z}^{(j)}(t)\right)^2 dt.
\]

We estimate the covariance function in Theorem 2 from empirical data. An estimate for
${\bf P}(X_{1i}\le s, \ X_{1j} \le t)$ is
\[
\frac{1}{n} \sum_{k=1}^n {\bf 1}\{\xi_{ki}\le \xi_{[ns],i}^{(i)}, \ \xi_{kj}\le \xi_{[nt],j}^{(j)}\}.
\]

It converges to the probability uniformely on $(s,t)$ in $[0,1]^2$. We estimate functions $L^{(j)}$ and 
matrix $G$ by their empirical counterparts.

Thus we construct a statistical test for accurate analysis of the data correspondence to the
linear regression model. This test allows one to use multiple
ordering of the initial multidimensional data and, due to this, to find non-obvious differences of the investigated
data from the model.

\vspace{20 mm}

{\bf Acknowledgement}

The work is supported by Mathematical Center in Akademgorodok 
under agreement No. 075-15-2019-1675 with the Ministry of Science and Higher Education of the Russian Federation.

\end{document}